\documentclass{amsart}

\usepackage{amscd,amssymb}

\def\Z{{\mathbb Z}}
\def\Q{{\mathbb Q}}
\def\R{{\mathbb R}}
\def\N{{\mathbb N}}
\def\C{{\mathbb C}}
\def\P{{\mathbb P}}
\def\V{{\mathbb V}}

\def\A{{\mathcal A}}
\def\H{{\mathcal H}}
\def\J{{\mathcal J}}
\def\M{{\mathcal M}}

\def\G{{\Gamma}}

\def\g{{\mathfrak g}}

\def\p{{\mathfrak p}}
\def\s{{\mathfrak s}}

\def\sp{\s\p}					
\def\kbar{\overline{k}}
\def\dec{\mathrm{dec}}
\def\Xred{X^\dec}
\def\Zred{Z^\dec}
\def\Ared{\A^\dec}
\def\Jred{\J^\dec}

\def\dot{{\bullet}}
\def\over{{\backslash}}

\def\Xo{\smash[t]{\overset{\circ}{X}}}
\def\Zo{\smash[t]{\overset{\circ}{Z}}}

\newcommand\im{\operatorname{im}}		

\newcommand\Jac{\operatorname{Jac}}
\newcommand\Pic{\operatorname{Pic}}

\newtheorem{theorem}{Theorem}[section]
\newtheorem{lemma}[theorem]{Lemma}
\newtheorem{proposition}[theorem]{Proposition}
\newtheorem{corollary}[theorem]{Corollary}
\newtheorem{bigtheorem}{Theorem}

\theoremstyle{definition}
\newtheorem{definition}[theorem]{Definition}

\theoremstyle{remark}
\newtheorem{remark}[theorem]{Remark}

\begin{document}

\title[Locally Symmetric Families of Curves and Jacobians]%
{Locally Symmetric Families of Curves and Jacobians}

\author{Richard Hain}

\address{Department of Mathematics\\ Duke University\\
Durham, NC 27708-0320}

\email{hain@math.duke.edu}

\thanks{This work was supported in part by a grant from the
National Science Foundation.}

\date{March 1997}

\maketitle

\section{Introduction}

In this paper we study locally symmetric families of curves and jacobians.
For the purposes of this paper, a jacobian is an abelian variety that is
$\Pic^0$ of a semi-stable curve. Denote the moduli space of principally
polarized abelian varieties of dimension $g$ by $\A_g$. It is a locally
symmetric variety. By a {\it locally symmetric family of jacobians}, we
mean a family of jacobians parameterized by a locally symmetric variety
$X$ where the period map $X \to \A_g$ is a map of locally symmetric
varieties  --- cf.\ (\ref{sym_jac}).  A {\it locally symmetric family
of curves} is a family of
semi-stable curves over a locally symmetric variety where $\Pic^0$ of each
curve in the family is an abelian variety and where the corresponding family
of jacobians is a locally symmetric family --- cf.\ (\ref{sym_curve}). Not
every locally symmetric family of jacobians can be lifted to a locally
symmetric family of curves, even if one passes to arbitrary finite
unramified covers of the base, as can be seen by looking at the universal
abelian variety of dimension 3 --- see (\ref{no_lift}).

\begin{definition}\label{bad}
A locally symmetric variety $X$ is {\it bad} if it has a locally symmetric
divisor. Otherwise, we shall say that $X$ is {\it good}.
\end{definition}

Each locally symmetric variety $X$ arises from a semisimple $\Q$-group
$G$ whose associated symmetric space is hermitian. We may suppose that $G$
is simply connected as a linear algebraic group. In this case we can write
$G$ as a product $G_1\times \dots \times G_m$ of almost simple $\Q$-groups.
There is a corresponding splitting $X' = X_1 \times \dots \times X_m$ of a
finite covering $X'$ of $X$, and $X$ is good if and only if each $X_i$ is
good. Since each $G_i$ is almost simple, each $X_i$ is also ``simple''
in the sense that it has no finite cover which splits as a product of
locally symmetric varieties.

\begin{remark}\label{no_divisors}
\begin{enumerate}
\item Suppose that $X$ is a locally symmetric variety whose associated
$\Q$-group $G$ is almost simple. Nolan Wallach \cite{wallach} has
communicated to me a proof that if the $\Q$-rank of $G$ is $\ge 3$,
then $X$ is good.
\item It follows from the classification of hermitian
symmetric spaces \cite[p.~518]{helgason} that if the non-compact factor
of $G(\R)$ is simple and not $\R$-isogenous to $SO(n,2)$ or $SU(n,1)$ for
any $n$, then $X$ is also good --- see (\ref{no_div}).
\end{enumerate}
\end{remark}

Recall that the symmetric space $SU(n,1)/U(n)$ is the complex $n$-ball.

\begin{bigtheorem}\label{curves}
Suppose that $p : C \to X$ is a non-constant locally symmetric family of
curves over the locally symmetric variety $X$. Suppose that the
corresponding $\Q$-group is almost simple. If every
fiber of $p$ is smooth, or if $X$ is is good and the generic
fiber of $p$ is smooth, then $X$ is a quotient of the complex $n$-ball.
\end{bigtheorem}

Note that the simplicity of the group $G$ implies that the period map
$X\to A_g$ is finite.

After studying the geometry of the locus of jacobians, we are able
to prove the following result about a locally symmetric family $B\to X$
of jacobians. In this result, $\Xred$ denotes the set of points of $X$
where the fiber is the jacobian is that of a singular curve
and $X^*$ denotes its complement. The locus of hyperelliptic jacobians
in $X$ will be denoted by $X_\H$, and $X^*\cap X_\H$ will be denoted by
$X_\H^*$.

\begin{bigtheorem}\label{jacobians}
Suppose that $B \to X$ is a non-constant locally symmetric family of
jacobians and that $X^*$ is non-empty. Suppose that the corresponding
$\Q$-group is almost simple. If $X$ is good, then either:
\begin{enumerate}
\item $X$ is a quotient of the complex $n$-ball, or
\item $g\ge 3$, each component of $\Xred$ is of complex codimension $\ge 2$,
$X_\H^*$ is a non-empty divisor in $X^*$ which is smooth if $g>3$, and the
family does not lift to a locally symmetric family of curves.
\end{enumerate}
\end{bigtheorem}

Of course, we can also apply this theorem to all locally symmetric 
subvarieties of $X$.
Note that the second case in the theorem does occur; take $g=3$ and
$X$ to be $\A_3$ (cf.\ (\ref{no_lift})). I do not know any examples where the
second case of the theorem holds and $g>3$, but suspect there are none.

These results depend upon the following kind of rigidity result for
mapping class groups. We shall denote the mapping class group associated
to a surface of genus $g$ with $n$ marked points and $r$ non-zero
boundary components by $\G_{g,r}^n$.

\begin{bigtheorem}\label{homom}
Suppose that $G$ is a simply connected, almost simple $\Q$-group,
that the symmetric space associated to $G(\R)$ is hermitian, and
that $\G$ is an arithmetic subgroup of $G$. If $G(\R)$ has real rank
$>1$ (that is, it is not $\R$-isogenous to the product of $SU(n,1)$
with a compact group) and if $g\ge 3$, then the
image of every homomorphism $\G \to \G_{g,r}^n$ is finite.
\end{bigtheorem}

This is a special case of a much more general result which was
conjectured by Ivanov and proved by Farb and Masur \cite{farb-masur}.
(They do not assume hermitian symmetric nor that $\G$ is arithmetic.) 
Their methods are very different from ours, so we have included a
proof of this theorem in the case when $r+n>0$, and a weaker statement
that is sufficient for our applications in the case $r+n=0$.

Our approach is via group cohomology. The basic idea is that if
$C \to X$ is a locally symmetric family of curves, then there
is a homomorphism from $\pi_1(X,*)$ to the mapping class group
$\G_g$, which is the orbifold fundamental group of the moduli
space of curves. The rigidity theorem above shows that there
are very few homomorphisms from arithmetic groups to mapping
class groups. To get similar results for locally symmetric families
of jacobians, we exploit the fact that the fundamental group of a
smooth variety does not change when a subvariety of codimension
$\ge 2$ is removed. It is for this reason that we are interested
in good locally symmetric varieties.

This work arose out of a question of Frans Oort, who asked if
there are any positive dimensional Shimura subvarieties of the jacobian
locus not contained in the locus of reducible jacobians. He suspects
that once the genus is sufficiently large, there may be none. His
interest stems from a conjecture of Coleman \cite[Conj.~6]{coleman}
which asserts that the set of points in $\M_g$ whose jacobians
have complex multiplication is finite provided that $g\ge 4$. This
is false as stated as was shown by de~Jong and Noot \cite{dejong} who
exhibited infinitely many smooth curves of genus 4 and genus 6 whose
jacobians have complex multiplication. (Note that when $g \le 3$ there
are infinitely many curves whose jacobian has complex multiplication
as $\M_g$ is dense in $\A_g$ in these cases.)
However, Oort tentatively believes Coleman's conjecture may be true
when $g$ is sufficiently large. Partial results on Oort's question,
which are complementary to results in this paper, have been obtained
by Ciliberto, van~der~Geer and Teixidor~i~Bigas in 
\cite{cili-vdg} and \cite{cili-vdg-tex}.

I am very greatful to Frans Oort for asking me this question, and
also to my colleagues Les Saper, Mark Stern (especially) and Jun Yang
for helpful discussions about symmetric spaces, arithmetic groups and
their cohomology. I would also like to thank Nolan Wallach for 
communicating his result to me, and Nicholai Ivanov for telling me
about his belief that rigidity for mapping class groups held in this
form. Finally, thanks to Ben Moonen, Carel Faber and Gerard van der
Geer for their very helpful comments on this manuscript, which helped
improve the exposition and saved me from many a careless slip.

\section{Background and Definitions}

Here we gather together some definitions to help the reader. A
{\it locally symmetric variety} is a locally symmetric space of the form
$$
\G \over G(\R) / K
$$
where $G$ is a semi-simple $\Q$-group, $\G$ is an arithmetic subgroup
of $G$, $K$ is a maximal compact subgroup of $G(\R)$, and the
symmetric space $G(\R)/K$ is hermitian. The moduli space of abelian
varieties with a level $l$ structure
$$
\A_g[l] = Sp_g(\Z)[l]\over Sp_g(\R)/U(g).
$$
is a primary example. Here $Sp_g(\Z)[l]$ denotes the level $l$ subgroup
of $Sp_g(\Z)$. This moduli space is well known to be a quasi-projective
variety. A fundamental theorem of Baily and Borel \cite{bb} asserts that
every locally symmetric variety is a quasi-projective algebraic variety.
The imbedding into projective space is given by automorphic forms.

Suppose that
$$
X_1 = \G_1 \over G_1(\R) / K_1 \text{ and } X_2 = \G_2 \over G_2(\R) / K_2
$$
are two locally symmetric varieties. A map of locally symmetric varieties
$X_1 \to X_2$ is a map induced by a homomorphism of $\Q$-algebraic groups
$G_1 \to G_2$.

Recall that a semi-simple $k$-group $G$ is {\it almost simple} if its
adjoint form is simple. An algebraic $k$-group is {\it absolutely (almost)
simple} if the corresponding group over $\kbar$ is (almost) simple.

Every simply connected $\Q$-group $G$ is the product of simply connected,
almost simple $\Q$-groups:
$$
G = \prod_{i=0}^n G_i
$$
This implies that if $\G$ is an arithmetic 
subgroup of $G$, then 
$$
\prod_{i=0}^n \G \cap G_i
$$
has finite index in $\G$. It follows that if $X$ is the locally symmetric
space associated to $\G$ and $X_i$ the locally symmetric space associated
to $\G_i$, then the natural mapping
$$
\prod_{i=0}^n X_i \to X
$$
is finite. If $X$ is a locally symmetric variety (i.e., $G(\R)/K$ is
hermitian), then each $X_i$ will also be a locally symmetric variety.

We shall call a locally symmetric variety associated to an arithmetic subgroup
of an almost simple $\Q$-group a {\it simple} locally symmetric variety.
They have no finite covers which are products of locally symmetric varieties.
Also, since the associated $\Q$-group is almost simple, every map from
a simple locally symmetric variety to a locally symmetric variety will be 
either constant or finite.

We shall make critical use of a vanishing theorem proved by Raghunathan
\cite{raghunathan}. Suppose that $\G$ is an arithmetic subgroup
of an almost simple $\Q$-group $G$. If $G(\R)$ has real rank $\ge 2$, then
$$
H^1(\G,V) = 0
$$
for all finite dimensional representations of $G(\R)$.

The classification of hermitian symmetric spaces \cite[p.~518]{helgason}
implies that if 
$$
X = \G \over G(\R) / K
$$
is a simple locally symmetric variety where $G(\R)$ has real rank 1, then
the associated symmetric space $G(\R)/K$ is the complex $n$-ball
$SU(n,1)/U(n)$. This rank condition, and the corresponding failure of
Raghunathan's Vanishing Theorem, is the principal reason we cannot
handle ball quotients.

Margulis's Rigidity Theorem \cite{margulis} has the following important
consequence for mappings between locally symmetric varieties. Suppose that
$X_1$ and $X_2$ are locally symmetric varieties as above and that $X_1$
is simple. If the real rank of $G_1(\R)$ is 2 or more, then for every
homomorphism $\phi : \G_1 \to \G_2$, there is a finite covering
$$
\begin{CD}
\G_1' \over G_1(\R)/K_1 @>>> \G_1 \over G_1(\R)/K_1 \cr
@| @| \cr
X_1' @>>> X_1 \cr
\end{CD}
$$
of locally symmetric varieties corresponding to a finite index subgroup
$\G_1'$ of $\G_1$, and a map of locally symmetric varieties
$X_1' \to X_2$ which induces the restriction
$$
\phi|_{\G_1'} : \G_1' \to \G_2
$$
of $\phi$ to $\G_1'$.

Every arithmetic group $\G$ has a torsion free subgroup of finite index.
When $\G$ is torsion free, the locally symmetric variety $\G\over G(\R)/K$
is a $K(\G,1)$. That is, it has fundamental group $\G$ and all of its higher
homotopy groups are trivial.

The level $l$ subgroup $\G_{g,r}^n[l]$ of the mapping class group $\G_{g,r}^n$
is the kernel of the natural homomorphism $\G_{g,r}^n \to Sp_g(\Z/l\Z)$. It
is torsion free when $l\ge 3$. In this case it is isomorphic to the fundamental
group of $\M_{g,r}^n[l]$, the moduli space of smooth projective curves of genus
$g$ over $\C$ with $n$ marked points, $r$ non-zero tangent vectors and a level
$l$ structure. (As is customary, shall omit $r$ and $n$ when they are zero.)
This isomorphism is unique up to an inner automorphism. We shall often fix a
level $l\ge 3$ to
guarantee that $\A_g[l]$ and $\M_{g,r}^n[l]$ are smooth.

\section{Homomorphisms to Mapping Class Groups}
\label{homoms}

Suppose that $\G$ is a discrete group and that and that
$\phi : \G \to \G_{g,r}^n$ is a homomorphism. Denote the composite
of $\phi$ with the natural homomorphism $\G_{g,r}^n \to Sp_g(\Z)$
by $\psi$. Denote the $k$th fundamental representation of $Sp_g(\R)$
by $V_k$ ($1\le k \le g$). Each of these can be viewed as a $\G$
module via $\psi$.

\begin{theorem}\label{coho}
If  $H^1(\Gamma,V_3)=0$ and $g \ge 3$, then
$\psi^\ast : H^2(Sp_g(\Z),\Q) \to H^2(\G,\Q)$ vanishes.
\end{theorem}

\begin{proof}
It suffices to prove the result when $r=n=0$.
I'll give a brief sketch. We can consider the group $H^1(\G_g,V_3)$.
One has the cup product map
$$
H^1(\G_g,V_3)^{\otimes 2} \to H^2(\G_g,\Lambda^2 V_3).
$$
Since $V_3$ has an $Sp_g$ invariant, skew symmetric bilinear form,
we have a map
$$
H^2(\G_g,\Lambda^2 V_3) \to H^2(\G_g,\R).
$$
We therefore have a map
$$
c: H^1(\G_g,V_3)^{\otimes 2} \to H^2(\G_g,\R).
$$
It follows from Dennis Johnson's results \cite{johnson} (cf.\ 
\cite[(5.2)]{hain:normal}) that the left hand group has dimension 1
when $g\ge 3$, and from \cite[\S\S 7--8]{hain:comp} that this map is
injective. It is also known that the natural map
$$
H^2(Sp_g(\Z),\R) \to H^2(\Gamma_g,\R)
$$
is an isomorphism for all $g\ge 3$ and that both groups are isomorphic
to $\R$ (cf.\ \cite{harer}, \cite{morita}, \cite{mumford}.)
The result now follows from the commutativity of the diagram.
$$
\begin{CD}
H^1(\G_g,V_3)^{\otimes 2} @>c>> H^2(\G_g,\R) \cr
@VVV @VVV \cr
H^1(\G,V_3)^{\otimes 2} @>c>> H^2(\G,\R)
\end{CD}
$$
\end{proof}

\begin{corollary}\label{nosplit}
For all $g\ge 3$ and all $l > 0$, the homomorphism
$\Gamma_g[l] \to Sp_g(\Z)[l]$ does not split.
\end{corollary}

\begin{proof}
By a result of Raghunathan \cite{raghunathan}, if $g\ge 2$, $\G$ is a finite
index subgroup of $Sp_g(\Z)$ and $V$ is a rational representation of
$Sp_g(\R)$, then $H^1(\G,V)$ vanishes. In particular, 
$$
H^1(Sp_g(\Z)[l],V_3) = 0
$$
when $g\ge 3$. If there were a splitting of the canonical homomorphism
$\G_g \to Sp_g(\Z)[l]$, then, by Theorem~\ref{coho}, the homomorphism
$$
H^2(Sp_g(\Z),\Q) \to H^2(Sp_g(\Z)[l],\Q)
$$
induced by the inclusion would be trivial. But when $g\ge 3$, Borel's
stability theorem \cite{borel} (see also \cite[3.1]{hain:torelli})
implies that this mapping is an isomorphism. It follows that no such
splitting can exist.
\end{proof}

\begin{remark}
There are cases when the cohomological obstruction above vanishes
but there is still no lift. An example is given in the appendix.
\end{remark}

Suppose that $X$ is a topological space with fundamental group $\pi$
and that $\V$ is the local system corresponding to the $\pi$ module
$V$. It is a standard fact that there is a natural map
$$
H^k(\pi,V) \to H^k(X,\V)
$$
which is an isomorphism when $k\le 1$ and injective when $k=2$.
We thus have the following corollary.

\begin{corollary}\label{top}
Suppose that $g\ge 3$, $l\ge 3$ and that $\phi : X \to \M_g[l]$ is a
continuous map from a topological space $X$ to $\M_g[l]$.
If $H^1(X,\V_3)=0$, then the map
$$
H^2(\A_g,\Q) \to H^2(X,\Q)
$$
induced by the composition $X \to \A_g$ of $\phi$ with the period map
vanishes. \qed
\end{corollary}

\section{Maps of Lattices to Mapping Class Groups}

In this section we apply the results of the previous section to
the case where $\G$ is the arithmetic group associated to a locally 
symmetric variety. We begin by recalling a result of Yang \cite{yang}.

Suppose that $G$ is a semisimple group over a number field
$k$. Yang \cite[(3.3)]{yang} defined an invariant $\ell(G/k)$ of $G$ as
follows. Fix a minimal parabolic $P$ of $G$. Let $T$ be the corresponding
maximal $k$-split torus. Denote the relative root system of $G$ by $\Phi$,
the base of $\Phi$ corresponding to $P$ by $\Delta$, and the set of positive
relative roots by $\Phi_+$. Set
$$
\rho_P = \frac{1}{2}\sum_{\alpha \in \Phi_+} \alpha.
$$
Then we can define
\begin{multline*}
\ell(G/k) = \max\big\{q \in \N : 2 \rho_P - \sum_{\alpha \in J} \alpha 
\text{ is strictly dominant} \cr
\text{ for all subsets $J$ of $\Phi_+$ of cardinality $< q$} \big\}.
\end{multline*}

Denote the Lie algebra of $G$ by $\g$ and its maximal compact
subgroup by $K$.
The importance of this invariant lies in the following result
of Yang \cite{yang}.

\begin{theorem}[Yang]\label{yang} 
If $G$ is a semisimple $\Q$-group and $\G$ an arithmetic subgroup of
$G$, then the natural map $H^m(\g,K) \to H^m(\G,\R)$
is injective whenever $m \le \ell(G/\Q)$. \qed
\end{theorem}

We will need the following consequence:

\begin{lemma}\label{comp}
Suppose that $\G$ is an arithmetic subgroup of an almost simple $\Q$-group
$G$. If $G\neq SL_2(\Q)$, then $H^2(\g,K) \to H^2(\G,\R)$ is injective.
\end{lemma}

\begin{proof}
When the $\Q$-rank of $G$ is zero, the corresponding locally symmetric
space is compact and the relative Lie algebra cohomology injects.
So we suppose that $G$ has positive $\Q$-rank. Note that $\ell(G/\Q)\ge 1$
for all such $G$.

Every almost simple $\Q$-group is of the form $R_{k/\Q}G$, where
$G$ is an absolutely almost simple $k$-group (\cite[p.~46]{tits}).
It is not difficult to see that if $G$ is a $k$-group, then
$$
\ell((R_{k/\Q} G)/\Q) \ge [k:\Q]\, \ell(G/k).
$$
So, by Yang's Theorem, it suffices to prove that for every absolutely
almost simple
group $G$, $\ell(G/k) \ge 2$. This follows from Tits' classification
of the decorated Dynkin diagrams of such groups \cite[pp~54--61]{tits}
as the Dynkin diagram of each these reduced root system is connected and
has $\ge 2$ vertices --- use \cite[\S10.2]{humph}.
\end{proof}

\begin{proof}[Proof of Theorem~\ref{homom}]\footnote{In the case
$r+n=0$ we shall prove that the image of $\G \to \G_g \to Sp_g(\Q)$
is finite. If $\G \to \G_g$ is the monodromy representation of a family
of jacobians, the finiteness of the image of this homomorphism implies
the isotriviality of the family. This is all we shall need in the
sequel.}
Since every simply connected almost simple $\Q$-group $G$ is of the
form $R_{k/\Q}G'$ where $k$ is a number field and $G'$ is a simply
connected, absolutely almost simple $k$-group, the classification of
Hermitian symmetric spaces \cite[p.~518]{helgason} implies that either
$G(\R)$ is isogenous to the product of $SU(n,1)$ and a compact group,
or has real rank $\ge 2$.

By the superrigidity theorem of Margulis \cite{margulis}, we know
that, after replacing $\G$ by a finite index subgroup if necessary,
there is a $\Q$-group homomorphism $G \to Sp_g(\Q)$ that induces
the homomorphism $\G \to \G_{g,r}^n \to Sp_g(\Z)$. It follows that the
representation $V_3$ of $\G$ is the restriction of a representation
of $G$. Since the real rank of $G$ is at least 2, it follows from
Raghunathan's result \cite{raghunathan} that $H^1(\Gamma,V_3)$
also vanishes. So, by (\ref{coho}), $H^2(Sp_g(\Z),\R) \to H^2(\G,\R)$
vanishes. But we have the commutative diagram
$$
\begin{CD}
H^2(Sp_g(\Z),\C) @>>> H^2(\G,\C) \cr
@| @AApA \cr
H^2(\sp_g(\R),U(g);\C) @>>> H^2(\g(\R),K;\C)\cr
@| @| \cr
H^2(X,\C) @>>> H^2(Y,\C)
\end{CD}
$$
where $X$ is the compact dual of Siegel space $Sp_g(\R)/U(g)$, $Y$ is
the compact dual of $G(\R)/K$, and the vertical maps are the standard
ones. Consider the homomorphism $G\to Sp_g(\R)$. If it is trivial, then
the image of $\Gamma$ in $Sp_g(\Z)$ is finite. (Remember that we passed
to a finite index subgroup earlier in the argument.) We will show that
it must be trivial.

Since $G(\R)$ is semisimple, it has a finite cover which is a product
of simple groups $G_i$. Since the symmetric space associated to $G(\R)$
is hermitian, the symmetric space associated to each $G_i$ is also
hermitian. Denote the compact dual of the symmetric space associated
to $G_i$ by $Y_i$. Then $Y$ is the product of the $Y_i$. Since each
$G_i$ is simple, the induced map $G_i \to Sp_g(\R)$ is either trivial
or has finite kernel. It follows that the corresponding maps $Y_i \to X$
of compact duals are either injective or trivial. If $G \to Sp_g(\R)$
is non-trivial, there is an $i$ such that $G_i \to Sp_g(\R)$ is injective. 
The corresponding map $Y_i \to X$ of compact duals is injective. Since
$X$ and $Y_i$ are K\"ahler, and since $H^2(Sp_g(\Z),\C)$ is one
dimensional, the map $H^2(X,\C)\to H^2(Y_i,\C)$ must be injective
as the K\"ahler form cannot vanish on $Y_i$. It follows that the bottom
map in the diagram is also injective.  Combining this with (\ref{comp}),
we see that the map $p$ is injective, which in turn implies that the
top map is injective. But this contradicts (\ref{coho}). So
$G\to Sp_g(\R)$ must be trivial. This is the statement we set out to
prove when $r+n = 0$.

Since $G \to Sp_g(\R)$ is trivial, we conclude that $\G$ maps to the Torelli
group $T_{g,r}^n$. But this group is residually torsion free nilpotent
(i.e., injects into its unipotent completion --- see \cite[(14.9)]{hain:torelli}) when $r+n > 0$. Since $H_1(\G,\Q)=0$, the
image of $\G$ in every unipotent group over $\Q$ is trivial. It follows
that the image of $\G$ in the unipotent completion of $T_{g,r}^n$, and
therefore in $T_{g,r}^n$, is trivial.\footnote{It would be interesting to
know if $T_g$ is residually torsion free nilpotent when $g$ is sufficiently
large.}
\end{proof}

\begin{corollary}\label{smooth}
If we have a non-constant family of non-singular curves of genus $\ge 3$
over a locally symmetric variety whose corresponding $\Q$-group $G$ is almost
simple, then the base is a quotient of the complex $n$-ball.
\end{corollary}

\begin{proof}
The base $X$ of the family is $X = \G \over G(\R)/K$
where $G$ is the associated $\Q$-group, $\G$ is an arithmetic subgroup of
$G$, and $K$ is a maximal compact subgroup of $G(\R)$. By passing to a finite
index subgroup if necessary, we may assume that $\G$ is torsion free. Since
the family is non trivial, the homomorphism $\G \to Sp_g(\Z)$ induced by
the period map of the family is not finite. It follows from
Theorem~\ref{coho} (see also the footnote to its proof) that $G(\R)$ must
have real rank 1. This implies that $G(\R)/K$ is the complex $n$-ball.
\end{proof}

\section{Locally Symmetric Hypersurfaces in Locally Symmetric 
Varieties}

In order to extend (\ref{smooth}) to locally symmetric families of
stable curves, we will need to know when a locally symmetric
variety has a locally symmetric hypersurface. Recall that every almost
simple $\Q$-group is of the form $R_{k/\Q}G$ where $k$ is a number
field and $G$ is an absolutely almost simple group over $k$.

\begin{proposition}
If $G' = R_{k/\Q}G$, where $G$ is a $k$-group, then $G(\R)$ is
hermitian symmetric space if and only if $k$ is a totally real
field and the symmetric space associated to each real imbedding
of $k$ into $\R$ is hermitian.
\end{proposition}

\begin{proof}
We have
$$
G'(\R) = \prod_{\nu : k \to \C} G_\nu.
$$
where $\nu$ ranges over complex conjugate pairs of imbeddings of
$k$ into $\C$. Note that $G_\nu$ is an absolutely simple real group
if $\nu$ is a real
imbedding and $G_\nu$ is not an absolutely simple real group if $\nu$
is not real.  The symmetric space associated to $G'(\R)$ is the product
of the symmetric spaces of the $G_\nu$. It is hermitian symmetric
if and only if each of its factors is. But if $\nu$ is not real,
then $G_\nu$ cannot be compact, and its symmetric space is not
hermitian (see the list in Helgason \cite[p.~518]{helgason}). So
$k$ is totally real.
\end{proof}

\begin{proposition}\label{no_div}
Suppose that $G$ is a simple real Lie group whose associated 
symmetric space is hermitian. If $G$ is not isogenous to $SO(n,2)$ or
$SU(n,1)$ for any $n$, then the symmetric space associated to $G$ has
no hermitian symmetric hypersurfaces.
\end{proposition}

\begin{proof}
This follows from the classification of irreducible hermitian symmetric
spaces \cite[p.~518]{helgason} using \cite[Theorem 10.5]{satake}.
\end{proof}

\section{Geometry of the Jacobian Locus}

The geometry of the locus of jacobians will play a significant role
in the sequel. The {\it jacobian locus} $\J_g[l]$ in $\A_g[l]$ is the
closure of the image of the period map $\M_g[l] \to \A_g[l]$. Denote the
locus of jacobians of singular curves by $\Jred_g[l]$, and the locus
$\J_g[l] - \Jred_g[l]$ of jacobians
of smooth curves by $\J_g^*[l]$. In this section we are concerned with
the geometry of $\J_g[l]$ along the hyperelliptic locus of $\J_g^*[l]$
and along $\Jred_g[l]$.

\begin{proposition}\label{cone}
If $g \ge 3$ and $l\ge 3$, then the locus of hyperelliptic jacobians in
$\J_g^*[l]$ is smooth, and its projective normal cone in $\A_g[l]$ at the
point $[C]$ is $\P(S^2(V_C))$, where $V_C$ is a vector space of dimension
$g-2$. Further, its projective normal cone in $\J_g^*[l]$ at the point $[C]$
is $\P(V_C)$, which is imbedded in $\P(S^2 V_C)$ via the ``Plucker imbedding.''
\end{proposition}

\begin{proof}
This follows directly from \cite{oort-steenbrink}.
\end{proof}

The following fact must be well known.

\begin{lemma}\label{no_lines}
If $V$ is a complex vector space, then there is no line $L$ in $\P(S^2V)$
that is contained in the image of the Veronese imbedding
$$
\nu : \P(V) \hookrightarrow \P(S^2V).
$$
\end{lemma}

\begin{proof}
Denote the hyperplane class of $\P(V)$ by $H_1$ and that of $\P(S^2V)$ by
$H_2$. Since a hyperplane section of $\P(V)$ in $\P(S^2V)$ is a quadric,
$\nu^\ast H_2 = 2 H_1$. If $L$ is a line in $\P(S^2V)$ that is contained
in $\P(V)$, then
$$
H_2|_L = 2H_1|_L \in H^2(L,\Z).
$$
This is impossible as $H_2|_L$ generates $H^2(L,\Z)$. So no such line can
exist.
\end{proof}

These two results combine to give us geometric information about how
smooth subvarieties of $\A_g[l]$ that are contained in $\J_g[l]$ intersect
the hyperelliptic locus.

\begin{proposition}\label{hyp}
Suppose that $g\ge 4$ and $l\ge 3$ and that $Z$ is a connected complex
submanifold of $\A_g[l]$. If $Z$ is contained in $\J_g^*[l]$, then one of
the following holds:
\begin{enumerate}
\item $Z$ is contained in the hyperelliptic locus,
\item $Z$ does not intersect the hyperelliptic locus,
\item the intersection of $Z$ with the hyperelliptic locus is a smooth
subvariety of $Z$ of pure codimension 1.
\end{enumerate}
\end{proposition}

\begin{proof}
Suppose that $Z$ is neither contained in the locus of hyperelliptic jacobians
nor disjoint from it. Suppose that $P$ is in the intersection of $Z$ and the
hyperelliptic locus. Denote by $V$ the $g-2$ dimensional vector space
corresponding to the point $P$ as in (\ref{cone}). Denote the quotient of
the tangent space $T_P Z$ by the Zariski tangent space of the intersection of
$Z$ with the hyperelliptic jacobians by $N_P$. Since $Z$ is contained in
$\J_g^*[l]$, it  follows from (\ref{cone}) that
$$
\P(N_P) \hookrightarrow \P(V) \stackrel{\nu}{\hookrightarrow} \P(S^2 V).
$$
Since $Z$ is a smooth subvariety of $\A_g[l]$, $\P(N_P)$ is a linear
subspace of $\P(S^2 V)$. It follows from (\ref{no_lines}) that $N_P$ has
dimension 0 or 1. If the dimension of $N_P$ is 0 for generic $P$, then
$Z$ is contained in the hyperelliptic locus. Since $Z$ is not contained
in the hyperelliptic locus, it must be that $N_P$ is one dimensional for
generic $P$ in the intersection. But since the dimension of $N_P$ is
bounded by 1, it is one for all $P$ in the intersection of $Z$ with the
locus of hyperelliptic jacobians. This implies that the intersection is
a smooth divisor.
\end{proof}

Denote the point in $\A_h[l]$ corresponding the the principally 
polarized abelian variety $A$ with level $l$ structure by $[A]$.
If $g' + g'' = g$, there is a generically finite-to-one morphism%
\footnote{This map is generically one-to-one if $g'\neq g''$ and 
generically two-to-one if $g'=g''$.}
$\mu: \A_{g'}[l] \times \A_{g''}[l] \to \A_g[l]$. It takes $([A'],[A''])$
to $[A'\times A'']$. Suppose that $[A'\times A'']$ is a simple point of the
image of $\mu$. Let
$$
V' = H^{0,1}(A') \text{ and } V'' = H^{0,1}(A'')
$$

\begin{lemma}
The projective normal cone to the image of $\im \mu$ at $[A\times B]$ is
naturally isomorphic to $\P(V'\otimes V'')$.
\end{lemma}

\begin{proof}
This follows directly by looking at period matrices.
\end{proof}

Now suppose that $C$ is a genus $g$ stable curve with one node and
whose normalization is $C' \cup C''$ where $g(C')=g'$ and $g(C'')=g''$.
Provided the level $l$ is sufficiently large, $[\Jac C]$ will be a
simple point of $\im \mu$. We now continue the discussion above with
$A' = \Jac  C'$ and $A'' = \Jac C''$. Note that we have canonical
imbeddings
$$
C' \to \P(V') \text{ and } C'' \to \P(V'').
$$
Composing the product of these with the Segre imbedding
$$
 \P(V') \times \P(V'') \to \P(V'\otimes V'')
$$
we obtain a morphism
\begin{equation}\label{directions}
\delta : C'\times C'' \to \P(V'\otimes V'')
\end{equation}
which is a finite onto its image when both $g'$ and $g''$ are 2 or
more, and is a projection onto the second factor when $g'=1$.

We can project the tangent cone of $\J_g[l]$ into the normal bundle of the
image of $\A_{g'}[l]\times \A_{g''}[l]$ in $\A_g[l]$. We shall call this image
the {\it relative normal cone} of $\J_g[l]$ in $\A_g[l]$. We can then
projectivize to get a map from the projectivized relative normal cone of
$\J_g[l]$ to the projective normal cone of $\A_{g'}[l]\times \A_{g''}[l]$ in
$\A_g[l]$.

\begin{proposition}
The image of the projectivized relative normal cone of $\J_g[l]$ in
$\A_g[l]$ at $\Jac(C'\times C'')$ in the projective normal cone of
$\A_{g'}[l]\times \A_{g''}[l]$ in $\A_g[l]$ at $\Jac(C'\times C'')$
is the image of the map (\ref{directions}).
\end{proposition}

\begin{proof}
This follows directly from \cite[Cor.~3.2]{fay}.
\end{proof}

Suppose that $Z$ is a complex submanifold of $\A_g[l]$ which is
contained in the jacobian locus. Set $\Zred = Z \cap \Jred[l]$. Then
$$
	\Zred = \bigcup_{j=1}^{[g/2]} \Zred_j
$$
where $\Zred_j$ is the intersection of $Z$ with the image of 
$\A_j[l]\times \A_{g-j}[l]$ in $\A_g[l]$. Set
$$
\Zo^\dec_j = \Zred_j - \bigcup_{i\neq j} \Zred_i.
$$

\begin{corollary}\label{cod_red}
Suppose that $g\ge 4$ and $1\le j \le g/2$ and that $Z$ is a complex
submanifold of $\A_g[l]$ that is contained in $\J_g^*[l]$. If $Z$
is not contained in $\Jred$ and $j\neq 2$, then $\Zo^\dec_j$ is empty or
is smooth and has pure codimension 1 in $Z$. If $\Zo^\dec_2$
is non-empty, it has codimension at most 2 in $Z$.
\end{corollary}

\begin{proof}
As in the proof above, it follows from the fact that the image of
$\delta$ contains no lines in $\P(V'\otimes V'')$ except when $j=2$.
This follows from the fact that canonical curves contain no lines
except in genus 2.
\end{proof}

\begin{remark}
One can find further restrictions on the largest strata of $\Zred$,
but so far I have not been able to find a use for them. Partial
results in this direction have been obtained in \cite{cili-vdg} and
\cite{cili-vdg-tex}.
\end{remark}

\section{Locally Symmetric Families of Curves}

\begin{definition}\label{sym_curve}
A {\it locally symmetric family of curves} is a family of stable
curves $p:C \to X$ where
\begin{enumerate}
\item $X$ is a locally symmetric variety;
\item the Picard group $\Pic^0$ of each curve in the family is an
abelian variety (i.e., the dual graph of each fiber is a tree);
\item the period map $X \to \A_g$ is a map of locally symmetric
varieties (i.e., it is induced by a $\Q$-algebraic group homomorphism).
\end{enumerate}
\end{definition}

We lose no generality by assuming that the base $X$ is a simple locally
symmetric variety.

By (\ref{smooth}), every locally symmetric family of curves
with all fibers smooth has a ball quotient as base. It
becomes more difficult to understand locally symmetric families when
the generic fiber is smooth, but some fibers are singular. Let
$\Xred_j$ be the closure of the locus in $X$ where the fiber has
two irreducible components, one of genus $j$, the other of genus
$g-j$. We will assume that $j \le g/2$. Set
$$
\Xo^\dec_j = \Xred_j - \bigcup_{i\neq j} \Xred_i.
$$

\begin{proposition}
Suppose that $X$ is simple and that $C \to X$ is a locally symmetric
family of curves of genus $g$ over whose generic fiber is smooth. If
$g \ge 4$, then each $\Xred_j$ is a locally symmetric subvariety of $X$.
Moreover, $\Xo^\dec_j$ is a locally symmetric hypersurface when $j\neq 2$,
and $\Xo^\dec_2$ has codimension at most 2.
\end{proposition}

\begin{proof}
Since $X$ is simple, the period map $X \to \A_g$ is constant or finite
onto its image. If it is constant, the result is immediate, so we 
assume that it is finite. Note that $\Xred$ is just the preimage
of $\Ared_g$, the locus of reducible principally polarized abelian
varieties. Since the period map $X \to \A_g$ is a map of  locally
symmetric varieties, each $\Xred_j$ is a locally symmetric subvariety
of $X$. The final statement follows from (\ref{cod_red}).
\end{proof}

\begin{proof}[Proof of Theorem~\ref{curves}]
The theorem is easily seen to be true when $g\le 2$ as all simple, locally
symmetric varieties of complex dimension $\le 3$ and real rank $\le 2$ are
ball quotients. Thus we can assume that $g\ge 3$. Fix a level $l\ge 3$
so that $\A_g[l]$ is smooth. By replacing $X$ by a finite covering,
we may assume that $X$ is smooth and that the family $C \to X$ is a family
of curves with a level $l$ structure, so that we have a period map
$X \to \A_g[l]$. The fundamental group $\G$ of $X$ is a torsion free
arithmetic group in a simple $\Q$-group $G$. Since the period map
is not constant, the homomorphism $\G \to Sp_g(\Z)$ induced by it
has infinite image. (If the image were finite, the family would be
isotrivial.)

If all fibers of $p : C \to X$ are smooth, then we have a lift of the
homomorphism $\G \to Sp_g(\Z)$ to a homomorphism $\G \to \G_g$. If $X$ is
good, $\Xred$ has codimension $\ge 2$, so that the inclusion of $X-\Xred$
into $X$ induces an isomorphism $\pi_1(X-\Xred,*) \cong \G$. Since we have a
map $X-\Xred \to \M_g[l]$, we have a homomorphism $\G \to \G_g$ which lifts
the homomorphism induced by the period map. Thus in both cases, we have
a lift $\G \to \G_g$ of the homomorphism induced by the period map.

Since the image of $\G$ in $Sp_g(\Z)$ is not finite,
Theorem~\ref{homom} implies that the corresponding real group $G(\R)$ is
isogenous to the product of $SU(n,1)$ with a compact group. It follows that
the symmetric space $G(\R)/K$ of $X$ is the complex $n$-ball.
\end{proof}

\section{Locally Symmetric Families of jacobians}

\begin{definition}
A {\it locally symmetric family of abelian varieties} is a family
$B \to X$ of principally polarized abelian varieties where the 
corresponding map $X \to \A_g$ is a map of locally symmetric varieties.
We will call such a family {\it essential} when the period map
is generically finite (and therefore finite) onto its image.
\end{definition}

The universal family of abelian varieties over $\A_g[l]$ is essential.
If $X$ is simple, then every locally symmetric family of abelian varieties
over $X$ is either essential or constant.

\begin{definition}\label{sym_jac}
A {\it locally symmetric family of jacobians} is a locally symmetric
family $B\to X$ of abelian varieties where the image of $B \to \A_g$
lies in the jacobian locus $\J_g$. It is {\it essential} when its 
family of jacobians is essential.
\end{definition}

The family of jacobians of a locally symmetric family of curves is
also a locally symmetric family of jacobians. It is natural to ask
whether every locally symmetric family of jacobians can be lifted
to a (necessarily locally symmetric) family of curves. The answer
is no. The problem is that the period map $\M_g[l]\to \A_g[l]$ is
$2:1$ and ramified along the hyperelliptic locus when $g\ge 3$.

\begin{proposition}\label{no_lift}
Suppose that $l\ge 3$. There is no family of curves over $\A_3[l]$ whose
family of jacobians is the universal family of abelian varieties over
$\A_3[l]$.
\end{proposition}

\begin{proof}
Let $U = \A_3[l] - \A_3^\dec[l]$. If the universal family of abelian
varieties over $\M_3[l]$ could be lifted to a family of curves, then one
would have a section of the period map $\M_3[l] \to U$. This is
impossible. One can explain this in two different ways:
\begin{enumerate}
\item When $l\ge 3$, the period map $\M_3[l] \to U$ is 2:1 and
ramified along the hyperelliptic locus. So it can have no section.
\item If there were such a section, one would have a splitting of
the group homomorphism $\G_3[l] \to \pi_1(U)$. But since each component
of $\A_3^\dec$ has codimension $\ge 2$ in $\A_3[l]$, it follows that
$\pi_1(U)$ is isomorphic to $\pi_1(\A_3[l])$, which is isomorphic to
$Sp_g(\Z)[l]$. The homomorphism $\G_3[l] \to \pi_1(U)$ does not split
by (\ref{nosplit}).
\end{enumerate}
\end{proof}

Suppose that $C \to X$ is an essential locally symmetric family of
jacobians (or curves). Set $X^* = X \cap \J_g^*[l]$.
Let $X_\H$ be the locus of points in $X$ where the fiber is 
hyperelliptic, and let $X_\H^* = X^* \cap X_\H$.

\begin{proof}[Proof of Theorem~\ref{jacobians}]
The theorem is easily seen to be true when $g \le 2$. So we supppose that
$g\ge 3$.
By passing to a finite index subgroup if necessary, we may suppose that
the arithmetic group $\G$ associated with $X$ is torsion free and that
the period map induces a homomorphism $\G \to Sp_g(\Z)[l]$ where $l\ge 3$.
With these assumptions, we have a period map $X\to \A_g[l]$ and both $X$ and
$\A_g[l]$ are smooth. Since $X$ is smooth and good, the inclusion
$X^* \hookrightarrow X$ induces an isomorphism on fundamental groups. Since
$\G$ is torsion free, $\pi_1(X,*)\cong \G$.

Now, if $X^*$ is contained in the locus of hyperelliptic jacobians, or
disjoint from it, then the period map $X \to \A_g[l]$ has a lift to
a map $X \to \M_g[l]$. That is, we can lift the family of jacobians over
$X^*$ to a family of curves. This implies that the homomorphism 
$\G \to Sp_g(\Z)$ induced by the period map, lifts to a homomorphism
$\G \to \G_g$. Since the period map is not constant, the image of $\G$ in
$Sp_g(\Z)$ is not finite (otherwise the family of abelian varieties would
be isotrivial). Applying Theorem~\ref{homom}, we see that $X$ is a quotient
of the complex $n$-ball.

The remaining case is where $X^*$ is neither disjoint from the locus 
of hyperelliptic jacobians nor contained in it. When $g=3$, the locus of
hyperelliptic curves is a divisor in $\M_g[l]$ from which it follows that
$X^*_\H$ is of pure codimension 1 in $X^*$. When $g > 3$, it follows
from (\ref{hyp}) that $X^*_\H$ is smooth of pure codimension 1 in $X^*$. Let
$Y$ be the fibered product
$$
Y = X \times_{\A_g[l]}\M_g[l]
$$
The map $Y \to X$ is 2:1 and ramified above the divisor $X^*_\H$ in $X^*$.
This implies that $Y \to X^*$ cannot be the restriction of a covering of
locally symmetric varieties as such coverings of $X$ are unramified as $\G$
is torsion free. 
\end{proof}

\begin{remark}
To try to eliminate the case where $X_\H$ is non-empty, one 
should study the second fundamental form (with respect to the
canonical metric of Siegel space) of the hyperelliptic locus
in $\A_g[l]$ and also that of its normal cone. If one is lucky,
this will give an upper bound on $\dim X_\H$, and
therefore  of $\dim X$ as well.
\end{remark}

\appendix
\section{An example}

Suppose that
$$
1 \to \Z \to G \to \G \to 1
$$
is a non-trivial central extension of a finite index subgroup $\G$ of
$Sp_g(\Z)$
by $\Z$. Suppose that the class of this extension is non-trivial in
$H^2(\G,\Q)$. One can ask whether the natural homomorphism $G \to Sp_g(\Z)$
lifts to a homomorphism $G \to \G_g$. We first show that the obstructions
of Theorem~\ref{coho} vanish when $g\ge 3$.

Since $g\ge 3$, $H^2(\G,\Q) =\Q$. Since the class of the extension is
non-trivial, an easy spectral sequence argument shows that $H^2(G,\Q)$
vanishes. Raghunathan's work \cite{raghunathan} implies that  that
$H^1(\G,V_3)$ vanishes where $V_3$ is the third fundamental representation
of $Sp_g(\R)$. Another easy spectral sequence argument implies that
$H^1(G,V_3)$ vanishes. The obstruction to a lift $G\to \G_g$ given in
Theorem~\ref{coho} therefore vanishes.

\begin{theorem}
There is an integer $g_0 \le 7$ such that
if $g\ge g_0$, there is no lift of the natural homomorphism
$p:G\to Sp_g(\Z)$ to a homomorphism $G \to \G_g$.
\end{theorem}

\begin{proof}
Suppose that the homomorphism $G \to \G$ lifts to a homomorphism
$\phi : G\to \G_g$. The latter homomorphism induces maps on cohomology
$$
H^\dot(\G_g,V) \to H^\dot(G,V)
$$
for all irreducible representations $V$ of $Sp_g(\R)$. Recall that the
cohomology group $H^1(\G_g,V_3)$ is isomorphic to $\R$. For each
$Sp_g(\R)$ equivariant map $p:V_3^{\otimes 6} \to \R$, there is a
natural map
$$
p_\ast : H^1(\G_g,V_3)^{\otimes 6} \to H^6(\G_g,\R).
$$
As is traditional, we denote the $i$th Chern class of the Hodge bundle
over $\A_g$ by $\lambda_i$. It follows from the work of Kawazumi and
Morita \cite{kw}, that for a suitable choice of $p$, the image of
$p_\ast$ is spanned by $\lambda_3$.

Since the diagram
$$
\begin{CD}
H^1(\G_g,V_3)^{\otimes 6} @>>> H^6(\G_g,\R) \cr
@VVV @VV{\phi^*}V \cr
H^1(G,V_3)^{\otimes 6} @>>> H^6(G,\R)
\end{CD}
$$
commutes, and since $H^1(G,V_3)$ vanishes when $g\ge 3$, it follows
that $\phi^*\lambda_3$ is trivial in $H^6(G,\R)$. We show that this
leads to a contradiction.

It follows from the work of Borel \cite{borel} that the 
ring homomorphism
$$
\R[\lambda_1,\lambda_3,\lambda_5,\dots] \to H^\dot(\G,\R)
$$
is an isomorphism in degrees $<g$. Since the characteristic class
of the extension $G$ of $\G$ by $\Z$ is a non-zero multiple of $\lambda_1$,
it follows that the ring homomorphism
$$
\R[\lambda_3,\lambda_5,\dots] \to H^\dot(G,\R)
$$
is an isomorphism in degrees $<g$. In particular, $\lambda_3$ is not
zero in $H^6(G,\R)$ when $g\ge 7$. It follows that there is no lifting
$\phi : G \to \G_g$ when $g\ge 7$.
\end{proof}

\end{document}